\documentclass[9pt,journal]{IEEEtran}
\usepackage{amsmath}
\usepackage{amsfonts}
\usepackage{latexsym}
\usepackage{amssymb}
\usepackage{cite}
\usepackage{array}
\usepackage{units}

\usepackage{ifpdf}
\ifCLASSINFOpdf
  \usepackage[pdftex]{graphicx}
  \graphicspath{{../figures/}}
  \DeclareGraphicsExtensions{.pdf,.jpeg,.png}
\else
  \usepackage[dvips]{graphicx}
  \graphicspath{{../figures/}}
  \DeclareGraphicsExtensions{.eps}
\fi

\newcommand{\norm}[1]{\left\|{#1}\right\|}

\newcommand{\mbf}[1]{\mathbf{#1}}

\newcommand{\LS}[0]{\mathrm{LS}}

\newcommand{\R}[0]{\mathcal{R}}
\newcommand{\C}[0]{\mathcal{C}}

\newcommand{\overbar}[1]{\mkern 1.5mu\overline{\mkern-1.5mu#1\mkern-1.5mu}\mkern 1.5mu}
\newcommand{\overtilde}[1]{\mkern 1.5mu\widetilde{\mkern-1.5mu#1\mkern-1.5mu}\mkern 1.5mu}

\begin{document}

\title{Comments on ``A Square-Root-Free Matrix Decomposition Method for Energy-Efficient Least Square Computation on Embedded Systems"}
\author{Mohammad~M.~Mansour,~\IEEEmembership{Senior~Member,~IEEE}
\thanks{M. M. Mansour is with the Department of Electrical and Computer Engineering at the American University of Beirut, Lebanon, e-mail: mmansour@ieee.org.}}

\maketitle

%
\begin{abstract}
A square-root-free matrix QR decomposition (QRD) scheme was rederived in~\cite{2014_Ren_sqrt_free} based on~\cite{1996_Bjork_numerical_methods} to simplify computations when solving least-squares (LS) problems on embedded systems. The scheme of~\cite{2014_Ren_sqrt_free} aims at eliminating both the square-root and division operations in the QRD normalization and backward substitution steps in the LS computations. It is claimed in~\cite{2014_Ren_sqrt_free} that the LS solution only requires finding the directions of the orthogonal basis of the matrix in question, regardless of the normalization of their Euclidean norms. MIMO detection problems have been named as potential applications that benefit from this. While this is true for unconstrained LS problems, we conversely show here that constrained LS problems such as MIMO detection still require computing the norms of the orthogonal basis to produce the correct result.
\end{abstract}
%

%
\section{Introduction}\label{sec_Introduction}
The problem of finding a vector $\mbf{x}\in\C^n$ such that $\mbf{y}\!=\!\mbf{A}\mbf{x}$ for a given \emph{matrix} $\mbf{A}\!\in\!\C^{m\times n}$ and \emph{observation vector} $\mbf{y}\!\in\!\C^m$ is well-studied (e.g., see~\cite{1996_Bjork_numerical_methods,1996_Golub}). For overdetermined systems ($m\geq n$), one minimizes the $2$-norm $\norm{\mbf{y}\!-\!\mbf{A}\mbf{x}}^2$ because it leads to tractable solutions; specifically, the distance quantity $\norm{\mbf{y}\!-\!\mbf{A}\mbf{x}}^2$ is differentiable, and the 2-norm is preserved under orthogonal transformations~\cite{1996_Golub}. Hence the least-``squares" problem finds $\hat{\mbf{x}}_{\LS}\!\in\!\C^n$ that gives
\begin{equation}\label{e:LS_problem}
  \underset{\mbf{x}\in \C^n}{\mathop{\min }} \norm{\mbf{y}\!-\!\mbf{A}\mbf{x}}^2.
\end{equation}

%
\section{Solving LS Problems with QRD}\label{sec:solving_LS_problems}
The \emph{full} QR decomposition of a matrix $\mbf{A}\!\in\!\C^{m\times n}$ is given by
\begin{equation}\label{e:QR_full}
  \mbf{A}_{m\times n} = \mbf{\overbar{Q}}_{m\times m}\mbf{\overbar{R}}_{m\times n},
\end{equation}
where $\mbf{\overbar{Q}}\!\in\!\C^{m\times m}$ is unitary (i.e., $\mbf{\overbar{Q}}^H\!\mbf{\overbar{Q}}\!=\!\mbf{I}_{m\times m}$, where $\mbf{I}_{m\times m}$ is the $m\!\times\!m$ identity matrix), and $\mbf{\overbar{R}}\!\in\!\C^{m\times n}$ is an upper-triangular matrix. The \emph{thin} QR decomposition of $\mbf{A}$ is given by
\begin{equation}\label{e:QR_thin}
  \mbf{A}_{m\times n} = \mbf{Q}_{m\times n}\mbf{R}_{n\times n},
\end{equation}
where $\mbf{Q}\!\in\!\C^{m\times n}$ has orthonormal columns, and $\mbf{R}\!\in\!\C^{n\times n}$ is a square upper-triangular matrix with real and positive diagonal entries. It is well-known~\cite{1996_Golub} that the two forms are related as follows:
\begin{align}\label{e:QR_thin_full_QR_thin}
  \mbf{A}
  \!=\! \mbf{\overbar{Q}}\mbf{\overbar{R}}
  \!=\! \left[
            \mbf{Q}_{m\times n}~~\mbf{\overtilde{Q}}_{m\times{(m-n)}}
        \right]
        \left[
            \begin{array}{c}
                \mbf{R}_{n\times n} \\
                \mbf{0}_{(m-n)\times n} \\
            \end{array}
        \right]
  \!=\! \mbf{Q}\mbf{R},
\end{align}
where $\mbf{\overtilde{Q}}_{m\times{(m-n)}}$ consists of the $(m\!-\!n)$ right-most columns of $\mbf{\overbar{Q}}$.

Assuming $\mbf{A}$ has full column rank, the unconstrained LS problem equivalently solves for $\mbf{x}$ by minimizing
\begin{align*}
  \underset{\mbf{x}\in \C^n}{\mathop{\min }} \norm{\mbf{y} \!-\! \mbf{A}\mbf{x}}^2
    \!&=\!
  \underset{\mbf{x}\in \C^n}{\mathop{\min }} \norm{\mbf{\overbar{Q}}^H\mbf{y}
    \!-\! \mbf{\overbar{Q}}^H\mbf{A}\mbf{x}}^2 \\
    \!&=\!
  \underset{\mbf{x}\in \C^n}{\mathop{\min }} \left\{ \norm{\mbf{Q}^H\mbf{y}- \mbf{R} \mbf{x}}^2 \right\} + \norm{\mbf{\overtilde{Q}}^H\mbf{y}}^2 \\
    \!&=\! \norm{\mbf{\overtilde{Q}}^H\mbf{y}}^2,
\end{align*}
Hence the minimizer of $\norm{\mbf{y} \!-\! \mbf{A}\mbf{x}}^2$ can be more simply obtained by solving the linear equations $\mbf{Q}^H\!\mbf{y}\!=\!\mbf{R}\mbf{x}$ using backward substitution
\begin{align*}
    x^{\star}\!=\!\underset{\mbf{x}\in \C^n}{\mathop{\arg\min }} \norm{\mbf{y} \!-\! \mbf{A}\mbf{x}}^2
    \!=\!
    \underset{\mbf{x}\in \C^n}{\mathop{\arg\min }} \norm{\mbf{Q}^H\mbf{y}- \mbf{R} \mbf{x}}^2
    \!=\! \mbf{R}^{-1}\mbf{Q}^H\mbf{y},
\end{align*}
which gives a zero cost function for $\underset{\mbf{x}\in \C^n}{\mathop{\min }} \left\{ \norm{\mbf{Q}^H\mbf{y} \!-\! \mbf{R} \mbf{x}}^2 \right\} \!=\!0$ and leaves a final residual of $\norm{\mbf{\overtilde{Q}}^H\mbf{y}}^2$. Note also that the thin form of the QRD suffices to obtain the unconstrained LS solution.

However, when solving \emph{constrained} LS problems over finite sets $\mathcal{X}^n \!\neq\! \mathcal{C}^n$, such as the case in MIMO detection problems over finite modulation constellations, one cannot guarantee that $x^{\star}\!=\!\mbf{R}^{-1}\mbf{Q}^H\mbf{y}$ actually belongs to $\mathcal{X}^n$, and the cost function $\underset{\mbf{x}\in \mathcal{X}^n}{\mathop{\min }}  \norm{\mbf{Q}^H\mbf{y} \!-\! \mbf{R} \mbf{x}}^2 $ is not necessarily 0 even if $\mbf{A}$ has full column rank. Many classes of the so-called ``hard-output'' MIMO detection algorithms find $\underset{\mbf{x}\in \mathcal{X}^n}{\mathop{\arg\min }}  \norm{\mbf{Q}^H\mbf{y} \!-\! \mbf{R} \mbf{x}}^2 $ (or an estimate of it) by actually computing the distance quantities $\norm{\mbf{Q}^H\mbf{y} \!-\! \mbf{R} \mbf{x}}^2$ and finding the minimum for all $\mbf{x}\!\in\!\mathcal{X}^n$ (when $n$ is small, e.g., 2 or 3), or for certain regions of points $\mbf{x}\!\in\!\mathcal{X}^n$ (e.g., see~\cite{2015_SP_mansour_256-QAM}), or over subsets of points in $\mathcal{X}^n$ that lie within some sphere centered around $\mbf{Q}^H\mbf{y}$.

Yet, in the more general ``soft-output" MIMO detection problems, the interest is in actually computing quantities of the form of \emph{differences of minimum distances}, known as log-likelihood ratios, over two disjoint partitions of the finite set $\mathcal{X}^n \!=\!  \mathcal{X}_1^n \dot{\cup} \mathcal{X}_2^n$:
\begin{multline*}
    \underset{\mbf{x}\in \mathcal{X}_1^n}{\mathop{\min }} \norm{\mbf{y} \!-\! \mbf{A}\mbf{x}}^2
    \!-\!
    \underset{\mbf{x}\in \mathcal{X}_2^n}{\mathop{\min }} \norm{\mbf{y} \!-\! \mbf{A}\mbf{x}}^2
    \!=\! \\
    \underset{\mbf{x}\in \mathcal{X}_1^n}{\mathop{\min }} \norm{\mbf{Q}^H\mbf{y}- \mbf{R} \mbf{x}}^2
    \!-\!
    \underset{\mbf{x}\in \mathcal{X}_2^n}{\mathop{\min }} \norm{\mbf{Q}^H\mbf{y}- \mbf{R} \mbf{x}}^2.
\end{multline*}

Here, the solution $x^{\star}\!=\!\mbf{R}^{-1}\mbf{Q}^H\mbf{y}$ of the original normal equations has no bearing on computing the above quantities. The application of the QRD for solving unconstrained LS problems in the context of MIMO detection is purely for complexity reduction reasons when enumerating points in $\mathcal{X}^n$ and computing distances. Furthermore, the thin form of the QRD suffices to find or estimate the minimum or arg-minimum distance quantities since the quantity $\norm{\mbf{\overtilde{Q}}^H\mbf{y}}$ originating from the full QRD is irrelevant. Hence we focus on this form of QRD in the discussion below.

%
\section{Comments and Discussion on~\cite{2014_Ren_sqrt_free}}\label{sec:discussion}
When performing the thin QRD using the modified Gram-Schmidt procedure~\cite{1996_Golub}, square-root operations are required to compute the diagonal entries of $\mbf{R}$ in~\eqref{e:QR_thin}, which represent the norms of the columns of $\mbf{A}$, as they are progressively normalized to become the columns of $\mbf{Q}$. It was first pointed out in~\cite{1996_Bjork_numerical_methods} and then elaborated further in~\cite{2014_Ren_sqrt_free}, that a square-root free QRD is possible by introducing a \emph{normalizer} diagonal matrix $\mbf{D}\!\in\!\R^{n\times n}$ into the so-called QDRD factorization as
\begin{equation}\label{e:QR_thin_expanded_with_D}
 \mbf{A} = \mbf{Q}\mbf{R} = \underbrace{\left(\mbf{Q}\mbf{D}^{-1}\right)}_{\triangleq \mbf{Q}'_{m\times n}} \underbrace{\mbf{D}^2\vphantom{\left(\mbf{Q}\mbf{D}^{-1}\right)}}_{\triangleq \mbf{D}'_{n\times n}} \underbrace{\left(\mbf{D}^{-1}\mbf{R}\right)}_{\triangleq \mbf{R}'_{n\times n}} = \mbf{Q}'\mbf{D}'\mbf{R}'.
\end{equation}
In~\eqref{e:QR_thin_expanded_with_D}, $\mbf{Q}'$ is an \emph{unnormalized} matrix with orthogonal columns since $\mbf{Q}'^H\mbf{Q}' \!=\! \mbf{D}^{-2} \!=\! \mbf{D}'^{-1} \!\neq\! \mbf{I}_{n\times n}$. Hence, $\mbf{Q}'\mbf{D}$ has orthonormal columns (unitary for $m\!=\!n$). Also, the matrix $\mbf{R}'$ is upper triangular, with unit diagonal elements.

It was shown in~\cite{2014_Ren_sqrt_free}, that to solve LS problems, one equivalently solves the modified normal equations as follows:
\begin{align*}
 \mbf{y} &= \mbf{Q'}\mbf{D}'\mbf{R}' \mbf{x}\\
 \left(\mbf{Q'D}\right)^H \mbf{y} &= \left(\mbf{Q'D}\right)^H  \mbf{Q'}\mbf{D}'\mbf{R}' \mbf{x} \\
 \mbf{D}\mbf{Q'}^H\mbf{y} &= \left(\mbf{Q'D}\right)^H  \left(\mbf{Q'}\mbf{D}\right) \mbf{D}\mbf{R}' \mbf{x}\\
 \mbf{D}\mbf{Q'}^H\mbf{y} &= \mbf{D}\mbf{R}' \mbf{x}\\
 \mbf{Q'}^H\mbf{y} &= \mbf{R}' \mbf{x}
\end{align*}
where in the last equation, $\mbf{D}$ ``cancels out", and hence is not needed in the solution.

As an immediate application to this, it was observed in~\cite[Section IV-D]{2014_Ren_sqrt_free} that the square-root-free QDRD can be applied to solve MIMO signal detection problems, which rely heavily on QR decompositions to reduce computational complexity (e.g., see~\cite{2015_SP_mansour_256-QAM}). This statement is true in the context of unconstrained LS problems over the $n$-dimensional complex field $\C^n$, where the LS solution is the solution of the normal equations as derived above according to~\cite{2014_Ren_sqrt_free}. We argue however that this observation is not accurate in the context of MIMO detection problems, and that the QDRD alone without the diagonal matrix $\mbf{D}$ is not adequate to solve constrained LS problems such as MIMO detection problems over finite sets $\mathcal{X}^n \!\neq\!\mathcal{C}^n$. In particular, the matrix $\mbf{D}$, which includes the square-root operations pertaining to the column norms that are to be eliminated by the scheme of~\cite{2014_Ren_sqrt_free}, are in fact still needed to find the constrained LS solution. This is due to the fact that MIMO detection problems require finding a \emph{constrained} LS solution in which the cost function $\underset{\mbf{x}\in \mathcal{X}^n}{\mathop{\min }}  \norm{\mbf{Q}^H\mbf{y} \!-\! \mbf{R} \mbf{x}}^2 $ is not necessarily 0, as discussed in Section~\ref{sec:solving_LS_problems}. The search space is constrained to be over a finite constellation of points and not over $\mathcal{C}^n$ where the solution of the normal equations applies.

Specifically, a MIMO system with $n$ transmit and $m\!\geq\!n$ receive antennas can be modeled by the equivalent complex baseband input-output system relation $\mbf{y}\!=\!\mbf{Ax}\!+\!\mbf{n}$, where $\mbf{A}\!\in\!\mathcal{C}^{m\times n}$ plays the role of a complex channel matrix, $\mbf{y}\!\in\!\mathcal{C}^{m}$ is the received complex signal vector, $\mbf{x}\!=\!\!\left[ x_1\ x_2\cdots x_n \right]^T\in\mathcal{X}^n$ is the $n\!\times\! 1$ transmitted complex symbol vector, and $\mbf{n}\!\in\!\mathcal{C}^{m}$ is a complex Gaussian circularly symmetric random noise vector. Each symbol $x_n$ belongs to a complex constellation $\mathcal{X}$. Finding the optimal solution (i.e., finding $\mbf{x}$ given $\mbf{y}$ and $\mbf{A}$) in the maximum-likelihood sense requires solving the following constrained LS problem
\begin{equation}\label{e:constrained_LS_problem}
  \underset{\mbf{x}\in \mathcal{X}^n}{\mathop{\arg\min }} \norm{\mbf{y}^{\vphantom{H}} - \mbf{A}\mbf{x}}^2 =
  \underset{\mbf{x}\in \mathcal{X}^n}{\mathop{\arg\min }} \norm{\mbf{Q}^H\mbf{y} - \mbf{R}\mbf{x}}^2
\end{equation}
over the finite $n$-dimensional constellation $\mathcal{X}^n$, where $\mbf{A}$ is factored according to~\eqref{e:QR_thin}. However, if the QDRD factorization in~\eqref{e:QR_thin_expanded_with_D} is applied, with the matrix $\mbf{Q'}$ being unnormalized (but still with orthogonal columns), then left-multiplying the quantity $\left(\mbf{y}\!-\!\mbf{A}\mbf{x}\right)$ by $\mbf{Q'}^H$ alone without $\mbf{D}$ does not preserve the ordering of the transformed squared-Euclidean distances $\norm{\mbf{Q}'^H\!\left(\mbf{y}\!-\!\mbf{A}\mbf{x}\right)}^2$ relative to the true $\norm{\mbf{Q}^H\!\left(\mbf{y}\!-\!\mbf{A}\mbf{x}\right)}^2$ distances, nor does it preserve the statistics of the noise vector $\mbf{n}$. In particular, we have
\begin{align}
  \underset{\mbf{x}\in \mathcal{X}^n}{\mathop{\arg\min }} \norm{\mbf{Q}^H\mbf{y} - \mbf{R}\mbf{x}}^2
  &= \underset{\mbf{x}\in \mathcal{X}^n}{\mathop{\arg\min }} \norm{\mbf{D}\mbf{Q'}^H\mbf{y} - \mbf{D}\mbf{R'}\mbf{x}}^2, \label{eq:distances_1}
\end{align}
but
\begin{align}
  \underset{\mbf{x}\in \mathcal{X}^n}{\mathop{\arg\min }} \norm{\mbf{D}\mbf{Q'}^H\mbf{y} - \mbf{D}\mbf{R'}\mbf{x}}^2
  &\neq \underset{\mbf{x}\in \mathcal{X}^n}{\mathop{\arg\min }} \norm{\mbf{Q'}^H\mbf{y} - \mbf{R'}\mbf{x}}^2. \label{eq:distances_2}
\end{align}
The last inequality follows from the fact that if
\begin{align*}
  \norm{\mbf{D}\mbf{Q'}^H\mbf{y} - \mbf{D}\mbf{R'}\mbf{x}_1}^2
  < \norm{\mbf{D}\mbf{Q'}^H\mbf{y} - \mbf{D}\mbf{R'}\mbf{x}_2}^2,
\end{align*}
for all vectors $\mbf{x}_2 \!\neq\! \mbf{x}_1 \! \in \mathcal{X}^n$,
then this does not necessarily imply that
\begin{align*}
  \norm{\mbf{Q'}^H\mbf{y} - \mbf{R'}\mbf{x}_1}^2
  < \norm{\mbf{Q'}^H\mbf{y} - \mbf{R'}\mbf{x}_2}^2.
\end{align*}
for all $\mbf{x}_2 \!\neq\! \mbf{x}_1 \! \in \mathcal{X}^n$. Viewed in terms of $\mbf{Q}'$ and $\mbf{R}'$, the left-hand side of~\eqref{eq:distances_2} is a \emph{weighted} constrained LS problem. Without proper weighting using $\mbf{D}$, the right-hand side of~\eqref{eq:distances_2} gives the wrong result. Hence the statement ``... \textit{which indicates that the normalization of factorization matrices in QRD is essentially redundant for solving LS problems}" in~\cite[pp. 74]{2014_Ren_sqrt_free} does not hold for constrained LS problems such as MIMO detection, and the scheme of~\cite{2014_Ren_sqrt_free} does not help in this case.

A further comment related to computations on energy-constrained embedded platforms, is that one might argue that although in~\eqref{eq:distances_1} the weighting of $\mbf{D}$ must be accounted for, the square-root operations themselves involved in $\mbf{D}$ are not needed because in the end the actual weights of the individual distance components of the vector 2-norm are the \emph{squares} of the diagonal entries of $\mbf{D}=[d_{ii}]$, i.e., $d_{ii}^2$ because squared-norm computations are involved, and not $d_{ii}$, which require square-root operations to be computed. However, this introduces $n$ extra multiplications per vector norm computation, times the number of vectors searched in $\mathcal{X}^n$, totaling an added complexity proportional to $n\cdot P$, where $P \gg n$ is the number of points searched in $\mathcal{X}^n$, while only $n$ square-root and $n$ division operations are saved due to applying the QDRD scheme. Without proper optimizations to handle these extra multiplications, this obviously defeats the whole purpose of eliminating the square-root and division operations via the QDRD scheme when used for example in well-known distance-based MIMO detection algorithms such as sphere decoding, K-best, sub-space detection, and LORD algorithms, among others.

%
\section{Conclusion}\label{sec:conclusion}
We have pointed out that the QDRD decomposition scheme that eliminates square-root operations to reduce computational complexity does not help in the context of constrained LS problems such as MIMO detection over finite sets. The normalization factors still play a role in the minimization process in this class of LS problems.

%
\bibliographystyle{IEEEtran}
\bibliography{IEEEabrv,sqrt_free_QR_Comment_bibfile}

\end{document}